\documentclass[11pt, a4paper]{article}
\usepackage[a4paper, margin=2.5cm]{geometry}
\usepackage[T1]{fontenc}
\usepackage[utf8]{inputenc}
\usepackage{indentfirst}
\usepackage{amsmath, amssymb, amsthm}
\usepackage{xcolor}
\definecolor{linknavy}{RGB}{22, 68, 128}
\usepackage[colorlinks=true, linkcolor=linknavy, citecolor=linknavy, urlcolor=linknavy]{hyperref}

\numberwithin{equation}{section}

\newtheorem{theorem}{Theorem}[section]
\newtheorem{lemma}[theorem]{Lemma}
\newtheorem{proposition}[theorem]{Proposition}
\newtheorem{corollary}[theorem]{Corollary}
\theoremstyle{definition}
\newtheorem{example}[theorem]{Example}
\theoremstyle{remark}
\newtheorem{remark}[theorem]{Remark}

\DeclareMathOperator{\Ric}{Ric}
\DeclareMathOperator{\Scal}{Scal}
\DeclareMathOperator{\Id}{Id}

\title{\textbf{A normal form for Randers metrics of scalar flag curvature}}
\author{Yalu Feng}
\date{}

\begin{document}

\maketitle

\begin{abstract}
In this paper, we prove that a smooth, strongly convex Randers metric in dimension $n \ge 3$ has scalar flag curvature if and only if it is locally the sum of an exact $1$-form and a Randers metric of scalar flag curvature with isotropic $S$-curvature, and we derive an explicit coordinate normal form for all such metrics. The proof proceeds by identifying the quotient Jacobi endomorphism of the magnetic spray with the null-screen curvature endomorphism of the canonical stationary Lorentzian lift, which shows that scalar flag curvature is equivalent to local conformal flatness of the lift and reduces the problem to timelike conformal Killing fields in Minkowski space. We also compute the Zermelo navigation data and characterize local isometries in terms of conformal transformations of Minkowski space.
\end{abstract}

\section{Introduction and main results}

A \emph{Randers metric} \cite{Randers1941} on a connected smooth $n$-manifold $M$ is a Finsler metric of the form
\[
F(x,y) = \alpha(x,y) + \beta_x(y), \qquad \alpha(x,y) = \sqrt{a_x(y,y)},
\]
where $a$ is a Riemannian metric and $\beta$ is a smooth $1$-form with $|\beta|_a < 1$. Its \emph{flag curvature} $K(x,y,\Sigma)$, defined for a nonzero flagpole $y \in T_x M$ and a $2$-plane $\Sigma \subset T_x M$ containing $y$, generalizes Riemannian sectional curvature; it is \emph{isotropic} if $K = K(x)$ depends only on $x$, and \emph{scalar} if $K = K(x,y)$ is independent of $\Sigma$. For $n = 2$, the flag curvature is automatically scalar, whereas isotropic flag curvature is not necessarily constant. For $n \ge 3$, Schur's theorem for Finsler metrics implies that
\[
\text{constant flag curvature} \iff \text{isotropic flag curvature} \implies \text{scalar flag curvature}.
\]
Unlike in Riemannian geometry---where, for a fixed plane $\Sigma$, $K$ is independent of the choice of flagpole $y$ in $\Sigma$ and all three conditions coincide for $n \ge 3$---the last implication is strict even for Randers metrics: $K(x,y)$ can still depend on $y$ at a fixed point, as in Example~\ref{ex:nonconst}.

The classification of Randers metrics with special flag curvature has been studied extensively through Zermelo navigation and projective Finsler geometry. Bao, Robles, and Shen \cite{BaoRoblesShen2004} classified strongly convex Randers metrics of \emph{constant} flag curvature in terms of Riemannian space forms and homothetic vector fields. Cheng and Shen \cite{ChengShen2009} extended this classification to metrics of scalar flag curvature with \emph{isotropic $S$-curvature}, showing that for $n \ge 3$ they arise via Zermelo navigation from a Riemannian metric of constant sectional curvature and a conformal vector field. Subsequently, Chen and Zhao \cite{ChenZhao2012} constructed an explicit family of non-projectively flat Randers metrics of scalar flag curvature with non-isotropic $S$-curvature, obtained by adding an exact $1$-form to a metric in the Cheng--Shen classification (see Example~\ref{ex:chen-zhao}). Motivated by such constructions, Shen and Xia \cite{ShenXia2013} coined the term \emph{projectively isotropic} for metrics formed by adding a closed $1$-form to a Randers metric of weakly isotropic flag curvature and provided further explicit examples (see Remark~\ref{rem:shen-xia-examples}). For general Randers metrics $F = \alpha + \beta$ without any $S$-curvature assumption, Shen and Civi \cite{ShenCivi2012} derived a Riemannian tensorial characterization of scalar flag curvature in terms of the curvature of $a$, the $2$-form $d\beta$, and its $a$-covariant derivative. In a complementary Lorentzian framework, Inagaki, Javaloyes, and Mar\c{c}al \cite{InagakiJavaloyesMarcal2026} recently characterized scalar flag curvature via the null Weyl curvature of an associated Lorentz--Finsler spacetime, yielding analogous curvature formulas for Randers metrics.

The purpose of this paper is to establish an explicit local normal form for smooth, strongly convex Randers metrics of scalar flag curvature in dimension $n \ge 3$, without any auxiliary assumptions.

\begin{theorem}\label{thm:main}
Let $F = \alpha + \beta$ be a smooth, strongly convex Randers metric on a smooth manifold $M$ of dimension $n \ge 3$. Then $F$ has scalar flag curvature if and only if each point admits a coordinate chart with image $U\subset\mathbb R^n$ in which there exist constants
$a_0, c, \sigma \in \mathbb{R}$, $b, p, q \in \mathbb{R}^n$,  $Q \in \mathfrak{so}(n)$, and a function $f \in C^\infty(U)$ such that, with respect to the standard Euclidean inner product $\langle \cdot, \cdot \rangle$, norm $|\cdot|$, and gradient $\nabla$ on $\mathbb{R}^n$, the fields
\begin{align}
v(x) &= a_0 + \langle b, x \rangle + c|x|^2, \label{eq:v} \\
w(x) &= p + Qx + \bigl(\sigma + 2\langle q, x \rangle\bigr) x - |x|^2 q, \label{eq:w} \\
D(x) &= v(x)^2 - |w(x)|^2 \label{eq:D}
\end{align}
satisfy $v > |w|$ and
\begin{equation}\label{eq:convexity}
2\langle w, \nabla f \rangle + v^2 |\nabla f|^2 - \langle w, \nabla f \rangle^2 < 1 \qquad \text{on } U,
\end{equation}
and $F$ is given by
\begin{equation}\label{eq:F}
F(x,y) = \frac{\sqrt{D(x)|y|^2 + \langle w(x), y \rangle^2} + \langle w(x), y \rangle}{D(x)} + df_x(y).
\end{equation}
\end{theorem}

The coordinate characterization in Theorem~\ref{thm:main} admits a clean geometric reformulation. To state it precisely, we fix the convention throughout the paper that the $S$-curvature of a Randers metric is computed with respect to its Busemann--Hausdorff volume density.

\begin{corollary}\label{cor:projectively-isotropic}
Let $F$ be a smooth, strongly convex Randers metric on a smooth manifold $M$ of dimension $n \ge 3$. Then $F$ has scalar flag curvature if and only if every point admits a neighborhood on which $F = F_0 + df$, where $f$ is a smooth function and $F_0$ is a strongly convex Randers metric with scalar flag curvature and isotropic $S$-curvature.
\end{corollary}

In the local coordinates of Theorem~\ref{thm:main}, the decomposition $F = F_0 + df$ is realized explicitly by taking $F_0$ to be the metric in \eqref{eq:F} with $f = 0$. In particular, $F$ is locally projectively isotropic in the sense of \cite{ShenXia2013}.

\vskip 3mm

The remainder of the paper is organized as follows. Section~\ref{sec:preliminaries} reviews the spray and magnetic formulations of scalar flag curvature and identifies the quotient Jacobi endomorphism with the null-screen curvature of the stationary Lorentzian lift $\hat{g} = a - (dt - \beta)^2$ (Lemma~\ref{lem:curvature-comparison}). Building on this identification, Section~\ref{sec:weyl} characterizes scalar flag curvature via the vanishing of the Weyl tensor of $\hat{g}$ (Lemma~\ref{lem:weyl-null-screen} and Proposition~\ref{prop:weyl}). Section~\ref{sec:normal-form} then combines the resulting local conformal flatness with an analysis of flat conformal Killing fields and a unit-time normalization along spacelike slices (Lemma~\ref{lem:unit-norm}) to prove Theorem~\ref{thm:main}. Addressing global aspects, Section~\ref{sec:global} characterizes local isometries via Minkowski conformal maps and establishes the global gluing description (Theorem~\ref{thm:local-iso} and Corollary~\ref{cor:global-atlas}). Finally, Section~\ref{sec:examples} provides explicit examples, derives the corresponding Zermelo navigation data, and proves Corollary~\ref{cor:projectively-isotropic}.

\section{Preliminaries}\label{sec:preliminaries}

We recall the standard spray formalism of Finsler geometry (see \cite{Shen2001} and \cite[Chapters~4--6]{Szilasi2010}). Let $\mathsf{L}= \frac{1}{2}F^2 = \frac{1}{2} F^2$ be the energy Lagrangian and $g_{ij}(x,y) = \mathsf{L}_{y^i y^j}(x,y)$ the fundamental tensor. The \emph{geodesic spray} of $F$ is the vector field
\[
\mathcal{S} = y^k \partial_{x^k} - 2G^k \partial_{y^k}  \qquad \text{on} \ TM \setminus \{0\},
\]
whose positively $2$-homogeneous coefficients $G^i(x,y)$ are determined by $2 g_{ij} G^j = \mathsf{L}_{y^i x^k} y^k - \mathsf{L}_{x^i}$ and govern geodesics via $\ddot{x}^i + 2G^i(x, \dot{x}) = 0$. Setting $N^i{}_j = G^i_{y^j}$, the \emph{Jacobi endomorphism} $R_y = R^i{}_j(x,y)\,\partial_{x^i} \otimes dx^j$ on $T_x M$ is defined by
\begin{equation}\label{eq:jacobi-endo}
R^i{}_j = 2 G^i_{x^j} - y^k G^i_{x^k y^j} + 2 G^k G^i_{y^k y^j} - G^i_{y^k} G^k_{y^j}.
\end{equation}
Along a geodesic $x(t)$ with velocity $y = \dot{x}(t)$, the variational field $J(t) \in T_{x(t)} M$ of a variation through affinely parametrized geodesics satisfies the linearized geodesic equation $\mathcal{D}^2 J + R_y J = 0$, where $(\mathcal{D}J)^i = \dot{J}^i + N^i{}_j J^j$ and the dot denotes differentiation with respect to $t$.

We collect two classical facts that will be used below.

\begin{lemma}[\cite{Shen2001, Szilasi2010}]\label{lem:quotient-scalar}
For every $y \in T_x M \setminus \{0\}$, the Jacobi endomorphism $R_y$ satisfies $R_y y = 0$ and is self-adjoint with respect to $g_y$. Consequently, $F$ has scalar flag curvature if and only if the endomorphism induced by $R_y$ on the quotient $T_x M / \mathbb{R}y$ is a scalar multiple of the identity for all $y \in T_x M \setminus \{0\}$.
\end{lemma}

\begin{lemma}[\cite{Shen2001, Szilasi2010}]\label{lem:proj-invariance}
If two positively $2$-homogeneous sprays are projectively related by $\tilde{G}^i = G^i + P(x,y)y^i$, where $P$ is positively $1$-homogeneous in $y$, then their Jacobi endomorphisms satisfy
\[
\tilde{R}_y \equiv R_y + (P^2 - \mathcal{S}P)\,\Id \pmod{\mathbb{R}y}.
\]
In particular, the induced endomorphisms on $T_x M / \mathbb{R}y$ differ by a scalar multiple of the identity.
\end{lemma}

\vskip 2mm

We next specialize to a Randers metric $F = \alpha + \beta$ and interpret its geodesics via magnetic flows and a stationary Lorentzian metric. For the magnetic and stationary Lorentzian viewpoints on Randers geometry, including the choice of spacelike time sections, see \cite[Sections~2.2--2.5]{GibbonsEtAl2009}. Let $\nabla^a$ denote the Levi-Civita connection of the Riemannian metric $a$, and define the $a$-skew-adjoint $(1,1)$-tensor field $A$ by
\begin{equation}\label{eq:A-def}
a(AX, Y) = \frac{1}{2} d\beta(X, Y) \qquad \text{for all } X, Y \in TM,
\end{equation}
where $d\beta(X, Y) = X(\beta(Y)) - Y(\beta(X)) - \beta([X, Y])$.
Let $\gamma_\varepsilon$ be a fixed-endpoint variation of $\gamma$, with velocity $Y=\dot\gamma$ and variation field $J$. The variation of the $\alpha$-length contributes $-\int_0^1 a(\nabla^a_Y(Y/\alpha),J)\,ds$. The variation of the line integral of $\beta$ contributes $\int_0^1 d\beta(J,Y)\,ds$. Thus
\[
\left.\frac{d}{d\varepsilon}\right|_{\varepsilon=0} \int_0^1 \bigl(\alpha(\dot{\gamma}_\varepsilon) + \beta(\dot{\gamma}_\varepsilon)\bigr)\,ds
= -\int_0^1 a\left(\nabla^a_{\dot{\gamma}} \frac{\dot{\gamma}}{\alpha(\dot{\gamma})} + 2A\dot{\gamma},\, J\right) ds,
\]
because $d\beta(J,Y)=-2a(AY,J)$. At constant $\alpha$-speed, the critical curves of the Randers length therefore satisfy the \emph{magnetic geodesic equation}
\begin{equation}\label{eq:mag-geo}
\nabla^a_{\dot{\gamma}} \dot{\gamma} = -2\alpha A\dot{\gamma}.
\end{equation}
Conversely, skew-adjointness of $A$ implies that every solution of \eqref{eq:mag-geo} has constant $\alpha$-speed. The corresponding spray has coefficients
\begin{equation}\label{eq:mag-spray}
G_{\mathrm{mag}} = G_a + \alpha Ay, \qquad y \in T_x M \setminus \{0\},
\end{equation}
where $G_a$ denotes the coefficients of the Riemannian geodesic spray. Thus, the magnetic spray and the geodesic spray of $F$ have the same oriented unparametrized geodesics. Because two sprays with the same oriented unparametrized geodesics differ only by a tangential term, their coefficients satisfy $G^i - G_{\mathrm{mag}}^i = P(x,y)y^i$, where positive $2$-homogeneity of $G$ and $G_{\mathrm{mag}}$ forces $P$ to be positively $1$-homogeneous in $y$.

\vskip 2mm

On the product manifold $\hat{M} = \mathbb{R} \times M$ with coordinate $t$ on the $\mathbb{R}$-factor, we introduce the \emph{canonical stationary Lorentzian lift} \cite{GibbonsEtAl2009}
\begin{equation}\label{eq:lift}
\hat{g} = a - (dt - \beta)^2, \qquad T = \partial_t.
\end{equation}
Since the coefficients of $\hat{g}$ are independent of $t$, the vector field $T$ is a unit timelike Killing vector field ($\hat{g}(T,T) = -1$). Moreover, the induced metric on the slice $\{t = 0\}$ is $a - \beta \otimes \beta$, which is positive-definite (so $\{t = 0\}$ is spacelike) precisely because $|\beta|_a < 1$. Since $\hat{g}(V, T) = -(dt - \beta)(V)$ for all $V \in T\hat{M}$, the $\hat{g}$-orthogonal distribution to $T$ is given by $T^\perp = \ker(dt - \beta)$. For any vector field $X \in \mathfrak{X}(M)$ (viewed on $\hat{M}$ via the product structure), its \emph{horizontal lift} $X^H = X + \beta(X)T$ satisfies $(dt - \beta)(X^H) = 0$ and hence lies in $T^\perp$; in fact, $T^\perp$ is spanned by such lifts. Moreover, because $(dt - \beta)(X^H) = 0$, the definition \eqref{eq:lift} immediately gives $\hat{g}(X^H, Y^H) = a(X, Y)$ for all $X, Y \in T_x M$, showing that $X \mapsto X^H$ defines an isometric isomorphism from $(T_x M, a)$ onto $(T_{(t,x)}^\perp, \hat{g})$. When no ambiguity arises, we omit the superscript $H$ in curvature expressions.

Recall that for any nonzero null vector $k \in T\hat{M}$ (i.e., $\hat{g}(k,k) = 0$), we have $\mathbb{R}k \subset k^\perp$. Because $\hat{g}$ has Lorentzian signature, its restriction to $k^\perp$ is positive semidefinite with radical $\mathbb{R}k$, and thus induces a positive-definite inner product on the quotient \emph{screen space} $k^\perp / \mathbb{R}k$. By the symmetries of the curvature tensor $\hat{R}$, the operator $X \mapsto \hat{R}(X, k)k$ satisfies $\hat{g}(\hat{R}(X, k)k, k) = 0$ and $\hat{R}(k, k)k = 0$; hence it maps $k^\perp$ into $k^\perp$ and maps $\mathbb{R}k$ to zero. It therefore induces a well-defined linear endomorphism on $k^\perp / \mathbb{R}k$, called the \emph{null-screen curvature endomorphism}.

We now compare this null-screen curvature endomorphism with the quotient Jacobi endomorphism of the magnetic spray \eqref{eq:mag-spray}.

\begin{lemma}\label{lem:curvature-comparison}
Let $u \in T_x M$ be an $a$-unit vector ($a(u,u) = 1$), let $\Pi_u : T_x M \to u^\perp$ denote the $a$-orthogonal projection, and set $k = T + u^H$. Under the canonical identifications
\[
T_x M / \mathbb{R}u \cong u^\perp \cong k^\perp / \mathbb{R}k,
\]
the quotient Jacobi endomorphism of the magnetic spray \eqref{eq:mag-spray} at $u$ coincides with the null-screen curvature endomorphism of $\hat{g}$ along $k$. Their common representation is the linear endomorphism $Q_u : u^\perp \to u^\perp$ given for $X \in u^\perp$ by
\begin{equation}\label{eq:Qu}
Q_u X = \Pi_u \Bigl( R^a(X, u)u - A^2 X + 3a(X, Au)Au + 2(\nabla^a_X A)u - (\nabla^a_u A)X \Bigr).
\end{equation}
\end{lemma}

\begin{proof}
We use the curvature convention $R(X,Y)Z = \nabla_X \nabla_Y Z - \nabla_Y \nabla_X Z - \nabla_{[X,Y]} Z$.

\emph{Step 1: Quotient Jacobi endomorphism of $G_{\mathrm{mag}}$.}
Consider a smooth variation of $\gamma$ through magnetic geodesics solving \eqref{eq:mag-geo}, with velocity field $Y = \dot{\gamma}$ and variation field $J$. Differentiating \eqref{eq:mag-geo} in the variation direction and using $\nabla^a_J Y = \nabla^a_Y J$ together with $J(\alpha) = a(Y/\alpha, \nabla^a_Y J)$, we obtain the Jacobi equation
\begin{equation}\label{eq:mag-jacobi}
\nabla^a_Y \nabla^a_Y J + R^a(J,Y)Y + 2\alpha A \nabla^a_Y J + 2a(Y/\alpha, \nabla^a_Y J)\, AY + 2\alpha (\nabla^a_J A)Y = 0.
\end{equation}
Set $\hat{u} = Y/\alpha$, $\hat{u}^\flat = a(\hat{u}, \cdot)$, and define the endomorphism $B = \alpha A + AY \otimes \hat{u}^\flat$. The dynamical covariant derivative along $\gamma$ associated with $G_{\mathrm{mag}}$ is $\mathcal{D} = \nabla^a_Y + B$, so that
\begin{align*}
\mathcal{D}J &= \nabla^a_Y J + BJ, \\
\mathcal{D}\mathcal{D}J &= \nabla^a_Y \nabla^a_Y J + 2B\nabla^a_Y J + (\nabla^a_Y B + B^2)J.
\end{align*}
Here $2B\nabla^a_Y J = 2\alpha A\nabla^a_Y J + 2a(\hat{u}, \nabla^a_Y J)AY$ matches the first-order terms in \eqref{eq:mag-jacobi}. Since $\alpha$ is constant along solutions of \eqref{eq:mag-geo}, we have $\nabla^a_Y \hat{u} = -2\alpha A\hat{u}$ and hence $\nabla^a_Y(AY) = (\nabla^a_Y A)Y - 2\alpha A^2 Y$. Using $a(\hat{u}, AY) = 0$ and $\hat{u}^\flat \circ A = -(A\hat{u})^\flat$, a direct computation gives
\begin{align*}
\nabla^a_Y B &= \alpha \nabla^a_Y A + \bigl((\nabla^a_Y A)Y - 2\alpha A^2 Y\bigr) \otimes \hat{u}^\flat - 2AY \otimes (AY)^\flat, \\
B^2 &= \alpha^2 A^2 + \alpha A^2 Y \otimes \hat{u}^\flat - AY \otimes (AY)^\flat,
\end{align*}
and therefore
\[
\nabla^a_Y B + B^2 = \alpha \nabla^a_Y A + \alpha^2 A^2 + \bigl((\nabla^a_Y A)Y - \alpha A^2 Y\bigr) \otimes \hat{u}^\flat - 3AY \otimes (AY)^\flat.
\]
Rewriting \eqref{eq:mag-jacobi} as $\mathcal{D}\mathcal{D}J + R_{\mathrm{mag}} J = 0$, we have
\[
R_{\mathrm{mag}} J = R^a(J,Y)Y + 2\alpha (\nabla^a_J A)Y - (\nabla^a_Y B + B^2)J.
\]
Since the right-hand side depends pointwise and linearly on $J$, $R_{\mathrm{mag}}$ defines a linear endomorphism of $T_x M$. Choosing the magnetic geodesic with initial velocity $Y = u$ at $x$ (so that $\alpha = \sqrt{a(u,u)} = 1$ and $\hat{u} = u$ at $x$), substituting the expression for $\nabla^a_Y B + B^2$, and applying $R_{\mathrm{mag}}$ to any tangent vector $X \in u^\perp$ (so that $u^\flat(X) = 0$), we obtain
\[
R_{\mathrm{mag}} X = R^a(X,u)u - A^2 X + 3a(X, Au)Au + 2(\nabla^a_X A)u - (\nabla^a_u A)X.
\]
Under the identification $T_x M / \mathbb{R}u \cong u^\perp$, the quotient Jacobi endomorphism at $u$ is represented by $\Pi_u \circ R_{\mathrm{mag}}|_{u^\perp}$, which gives \eqref{eq:Qu}.

\emph{Step 2: Null-screen curvature of $\hat{g}$.}
For $X, Y \in \mathfrak{X}(M)$, set $\varpi(X,Y) = a(AX, Y) = \frac{1}{2}d\beta(X,Y)$. Since $T = \partial_t$ commutes with vector fields lifted from $M$, we have $[X^H, T] = 0$ and
\[
[X^H, Y^H] = [X,Y] + \bigl(X(\beta(Y)) - Y(\beta(X))\bigr)T = [X,Y]^H + 2\varpi(X,Y)T.
\]
Applying the semi-Riemannian Koszul formula \cite[Chapter~3]{ONeill1983} together with $\hat{g}(X^H, Y^H) = a(X,Y)$, $\hat{g}(X^H, T) = 0$, and $\hat{g}(T,T) = -1$ yields the Levi-Civita connection $\hat{\nabla}$ of $\hat{g}$:
\[
\hat{\nabla}_{X^H} Y^H = (\nabla^a_X Y)^H + \varpi(X,Y)T, \qquad \hat{\nabla}_{X^H} T = \hat{\nabla}_T X^H = (AX)^H, \qquad \hat{\nabla}_T T = 0.
\]
Identify horizontal vectors with their projections to $TM$ and use $(\cdot)_{\mathrm{hor}}$ for the horizontal component under this identification. The connection formulas give
\begin{align*}
(\hat{R}(X^H,Y^H)Z^H)_{\mathrm{hor}} &= R^a(X,Y)Z + \varpi(Y,Z)AX - \varpi(X,Z)AY - 2\varpi(X,Y)AZ, \\
(\hat{R}(X^H,T)T)_{\mathrm{hor}} &= -A^2 X, \\
(\hat{R}(X^H,T)Y^H)_{\mathrm{hor}} &= (\nabla^a_X A)Y, \\
(\hat{R}(X^H,Y^H)T)_{\mathrm{hor}} &= (\nabla^a_X A)Y - (\nabla^a_Y A)X.
\end{align*}
For $k=T+u^H$, with $a(u,u)=1$ and $X\perp u$, expand $\hat R(X^H,k)k$ by bilinearity. Using $\varpi(X,u)=-a(X,Au)$ in the horizontal component gives
\begin{align*}
(\hat R(X^H,k)k)_{\mathrm{hor}}
 &= R^a(X,u)u-3\varpi(X,u)Au-A^2X+(\nabla^a_X A)u
       +\bigl((\nabla^a_X A)u-(\nabla^a_u A)X\bigr)\\
 &= R^a(X,u)u-A^2X+3a(X,Au)Au+2(\nabla^a_X A)u-(\nabla^a_u A)X .
\end{align*}
Every $v\in k^\perp$ decomposes uniquely as $v=ck+X^H$, with $c=-\hat g(v,T)$ and $X\in u^\perp$. Thus $X\mapsto[X^H]$ is an isometry from $u^\perp$ onto the screen. In this representation the screen curvature sends $X$ to $\Pi_u((\hat R(X^H,k)k)_{\mathrm{hor}})$, giving \eqref{eq:Qu}.
\end{proof}

At unit speed, with magnetic form $d\beta$ and Lorentz force $-2A$, the operator $Q_u$ in \eqref{eq:Qu} coincides with the magnetic curvature operator of \cite{Assenza2024,AssenzaReberTerek2025}, and its quadratic form $a(Q_u X, X)$ recovers the stationary curvature formula in \cite[Proposition~A.6]{InagakiJavaloyesMarcal2026}. Lemma~\ref{lem:curvature-comparison} identifies the quotient Jacobi endomorphism of the magnetic spray on $T_x M / \mathbb{R}u$ with the null-screen curvature endomorphism of the stationary lift on $k^\perp / \mathbb{R}k$, showing that both are represented by $Q_u$.

\section{Algebraic characterization via the Weyl tensor}\label{sec:weyl}

Let $g$ be a Lorentzian metric of dimension $m = n + 1 \ge 4$. Its Schouten tensor is
\[
\mathsf{P} = \frac{1}{m-2}\left(\Ric - \frac{\Scal}{2(m-1)} g\right),
\]
and its Riemann curvature tensor admits the algebraic Ricci--Weyl decomposition (see \cite[Chapter~7]{Lee2018} and \cite{ONeill1983})
\begin{equation}\label{eq:weyl-decomp}
R(X,Y)Z = \mathcal{W}(X,Y)Z + \mathsf{P}(Y,Z)X - \mathsf{P}(X,Z)Y + g(Y,Z)\mathsf{P}^\# X - g(X,Z)\mathsf{P}^\# Y,
\end{equation}
where $\mathcal{W}$ is the Weyl tensor, which has the same algebraic symmetries as $R$ and is completely trace-free. In particular, for any null vector $k$ and any $X, Y \in k^\perp$, the Schouten terms in \eqref{eq:weyl-decomp} contribute $\mathsf{P}(k,k)g(X,Y)$ to $g(R(X,k)k, Y)$.

\begin{lemma}\label{lem:weyl-null-screen}
Fix a point $p$ of a Lorentzian manifold $(M,g)$ of dimension $m\ge4$. For each nonzero null vector $k\in T_pM$, let
\[
\mathcal{R}_k \colon k^\perp / \mathbb{R}k \to k^\perp / \mathbb{R}k, \qquad [X] \mapsto [R(X,k)k],
\]
denote the null-screen curvature endomorphism along $k$. Then $\mathcal{W}_p=0$ if and only if $\mathcal R_k$ is scalar for every nonzero null $k\in T_pM$. In this case $\mathcal R_k=\mathsf{P}(k,k)\,\Id$.
\end{lemma}

\begin{proof}
For any nonzero null vector $k \in T_p M$ and any $X, Y \in k^\perp$, contracting \eqref{eq:weyl-decomp} on $(X,k,k,Y)$ yields
\begin{equation}\label{eq:screen-ricci-weyl}
g(R(X,k)k, Y) = \mathcal{W}(X,k,k,Y) + \mathsf{P}(k,k)g(X,Y),
\end{equation}
where $\mathcal{W}(X,Y,Z,U) := g(\mathcal{W}(X,Y)Z, U)$. Recall that $g$ descends to a positive-definite inner product $\bar{g}([X],[Y]) := g(X,Y)$ on the screen space $k^\perp / \mathbb{R}k$, with respect to which $\bar{g}(\mathcal{R}_k[X], [Y]) = g(R(X,k)k, Y)$.

\emph{Sufficiency.} If $\mathcal{W}_p = 0$, then \eqref{eq:screen-ricci-weyl} gives $\bar{g}(\mathcal{R}_k[X], [Y]) = \mathsf{P}(k,k)\bar{g}([X],[Y])$ for all $[X],[Y] \in k^\perp / \mathbb{R}k$. Nondegeneracy of $\bar{g}$ implies $\mathcal{R}_k = \mathsf{P}(k,k)\,\Id$ for every nonzero null vector $k \in T_p M$.

\emph{Necessity.} Suppose that $\mathcal{R}_k = \lambda_k\,\Id$ for every nonzero null vector $k \in T_p M$. By \eqref{eq:screen-ricci-weyl}, we have
\[
\mathcal{W}(X,k,k,Y) = \bigl(\lambda_k - \mathsf{P}(k,k)\bigr)g(X,Y) \qquad \text{for all } X, Y \in k^\perp.
\]
Choose a complementary spacelike subspace $\mathcal{E} \subset k^\perp$ such that $k^\perp = \mathbb{R}k \oplus \mathcal{E}$, and let $\{e_1, \dots, e_{m-2}\}$ be a $g$-orthonormal basis of $\mathcal{E}$. Since $\mathcal{E}^\perp$ is a $2$-dimensional Lorentzian subspace containing $k$, there exists a unique null vector $\ell \in \mathcal{E}^\perp$ satisfying $g(k,\ell) = -1$, yielding a null frame $\{k, \ell, e_1, \dots, e_{m-2}\}$ of $T_p M$. Because $\mathcal{W}$ is completely trace-free, taking the metric trace of $U \mapsto \mathcal{W}(U,k)k$ gives
\[
0 = -\mathcal{W}(\ell,k,k,k) - \mathcal{W}(k,k,k,\ell) + \sum_{i=1}^{m-2} \mathcal{W}(e_i,k,k,e_i) = (m-2)\bigl(\lambda_k - \mathsf{P}(k,k)\bigr),
\]
where the first two terms vanish by the skew-symmetry of $\mathcal{W}$. Since $m \ge 4$, it follows that $\lambda_k = \mathsf{P}(k,k)$, and hence
\begin{equation}\label{eq:weyl-null}
\mathcal{W}(X,k,k,Y) = 0 \qquad \text{for all } X, Y \in k^\perp.
\end{equation}

Now fix a unit timelike vector $T \in T_p M$ ($g(T,T) = -1$) and set $E := T^\perp$, which is a Euclidean inner product space of dimension $n = m - 1 \ge 3$ under $g|_E$. Define the symmetric bilinear form $H(X,Y) := \mathcal{W}(X,T,T,Y)$ on $E$. For any unit vector $u \in E$, the vectors $k_\pm := T \pm u$ are nonzero and null, and any $X, Y \in E \cap u^\perp$ lie in $k_+^\perp \cap k_-^\perp$. Expanding \eqref{eq:weyl-null} for $k = T + u$ and $k = T - u$ yields
\[
\mathcal{W}(X,T,T,Y) + \mathcal{W}(X,u,u,Y) \pm \bigl(\mathcal{W}(X,T,u,Y) + \mathcal{W}(X,u,T,Y)\bigr) = 0.
\]
Adding and subtracting these two identities separates the terms with even and odd numbers of $T$-arguments:
\begin{itemize}
\item \emph{Components with zero or two $T$-arguments.} Adding the two equations gives $\mathcal{W}(X,u,u,Y) = -H(X,Y)$ for all unit vectors $u \in E$ and all $X, Y \in E \cap u^\perp$. In particular, for any orthogonal pair of unit vectors $u, v \in E$, the pair symmetry of $\mathcal{W}$ implies
\[
H(u,u) = -\mathcal{W}(u,v,v,u) = -\mathcal{W}(v,u,u,v) = H(v,v).
\]
Since $\dim E \ge 3$, any two unit vectors in $E$ admit a common orthogonal unit vector; hence $H(u,u)$ is constant on the unit sphere of $E$, so $H(X,Y) = h\,g(X,Y)$ on $E$ for some $h \in \mathbb{R}$. Because $\mathcal{W}$ is trace-free and $\mathcal{W}(T,T,T,T) = 0$, taking the metric trace of $\mathcal{W}(\cdot, T, T, \cdot)$ gives $nh = 0$, so $H = 0$. Consequently, $\mathcal{W}(v,u,u,v) = 0$ for all orthogonal unit vectors $u, v \in E$, and hence for all $u, v \in E$ by bilinearity and skew-symmetry. Standard polarization together with the algebraic curvature symmetries and the first Bianchi identity of $\mathcal{W}$ then implies $\mathcal{W}|_{E^4} = 0$.

\item \emph{Components with one $T$-argument.} Subtracting the two equations and setting $X = Y = v \in E \cap u^\perp$ with $|u|_g = 1$ gives $\mathcal{W}(v,T,u,v) + \mathcal{W}(v,u,T,v) = 0$. By the pair symmetry and skew-symmetry of $W$, we have $\mathcal{W}(v,u,T,v) = \mathcal{W}(T,v,v,u) = \mathcal{W}(v,T,u,v)$, so $\mathcal{W}(v,T,u,v) = 0$ whenever $u \perp v$ (by homogeneity in $u$). For an arbitrary $u \in E$ and $v \in E \setminus \{0\}$, decomposing $u = u_\perp + c v$ with $u_\perp \perp v$ and using $\mathcal{W}(v,T,v,v) = 0$ yields $\mathcal{W}(v,T,u,v) = 0$ for all $u, v \in E$. Polarizing in $v$ (replacing $v$ by $v + z$) shows that the trilinear form $B(v,u,z) := \mathcal{W}(v,T,u,z)$ on $E$ is skew-symmetric in $v$ and $z$. Since $W$ is also skew-symmetric in its last two slots, $B$ is an alternating $3$-form on $E$. Applying the first Bianchi identity to $\mathcal{W}$ gives
\[
0 = \mathcal{W}(v,T,u,z) + \mathcal{W}(T,u,v,z) + \mathcal{W}(u,v,T,z) = 3B(v,u,z),
\]
whence $B = 0$.
\end{itemize}
Finally, with respect to the orthogonal decomposition $T_p M = \mathbb{R}T \oplus E$, multilinearity and skew-symmetry reduce every component of $\mathcal{W}_p$ to one with at most two $T$-arguments—namely, up to sign, a component of $\mathcal{W}|_{E^4}$, $B$, or $H$. Since all three vanish, we conclude that $\mathcal{W}_p = 0$.
\end{proof}

Combining Lemmas~\ref{lem:quotient-scalar}--\ref{lem:curvature-comparison} with Lemma~\ref{lem:weyl-null-screen} yields the following Lorentzian characterization of scalar flag curvature.

\begin{proposition}\label{prop:weyl}
Let $(M,F)$ be a smooth, strongly convex Randers manifold of dimension $n \ge 3$, where $F = \alpha + \beta$. Then $F$ has scalar flag curvature if and only if the Weyl tensor $\mathcal{W}$ of its canonical Lorentzian lift $(\hat{M}, \hat{g})$ in \eqref{eq:lift} vanishes identically.
\end{proposition}

\begin{proof}
Fix a point $x \in M$. Since $\hat{g}$ in \eqref{eq:lift} is invariant under translations in $t$, the Weyl tensor $\mathcal{W}$ vanishes identically on $\mathbb{R} \times \{x\}$ if and only if $\mathcal{W}_{(0,x)} = 0$.

By Lemma~\ref{lem:quotient-scalar} and positive $2$-homogeneity of the Jacobi endomorphism $R_y$ in $y$, $F$ has scalar flag curvature at $x$ if and only if the endomorphism induced by $R_y$ on $T_x M / \mathbb{R}y$ is a scalar multiple of $\Id$ for every $y \in T_x M \setminus \{0\}$, which by homogeneity is equivalent to the same condition for every $a$-unit vector $u \in T_x M$. Since the geodesic spray $G$ of $F$ and the magnetic spray $G_{\mathrm{mag}}$ in \eqref{eq:mag-spray} are projectively related, Lemma~\ref{lem:proj-invariance} implies that at each $a$-unit vector $u \in T_x M$, the quotient Jacobi endomorphism of $F$ on $T_x M / \mathbb{R}u$ is a scalar multiple of $\Id$ if and only if the same holds for the quotient Jacobi endomorphism of $G_{\mathrm{mag}}$. By Lemma~\ref{lem:curvature-comparison}, this is in turn equivalent to the null-screen curvature endomorphism $\mathcal{R}_k$ of $\hat{g}$ at $(0,x)$ being a scalar multiple of $\Id$ for every null vector of the form $k = T + u^H$ with $a(u,u) = 1$.

Now let $k \in T_{(0,x)}\hat{M}$ be an arbitrary nonzero null vector. Under the orthogonal decomposition $T_{(0,x)}\hat{M} = \mathbb{R}T \oplus T^\perp$, write $k = cT + v^H$ with $c \in \mathbb{R}$ and $v \in T_x M$. The null condition $0 = \hat{g}(k,k) = -c^2 + a(v,v)$ together with $k \neq 0$ forces $c \neq 0$; setting $u := c^{-1}v$ gives $a(u,u) = 1$ and $k = c(T + u^H)$. Because the curvature tensor $\hat{R}$ is multilinear, $\mathcal{R}_{ck} = c^2 \mathcal{R}_k$, so $\mathcal{R}_k$ is a scalar multiple of $\Id$ for every nonzero null vector $k \in T_{(0,x)}\hat{M}$ if and only if this holds for all $k = T + u^H$ with $a(u,u) = 1$. Applying Lemma~\ref{lem:weyl-null-screen} to $(\hat{M}, \hat{g})$ at $(0,x)$ therefore shows that $F$ has scalar flag curvature at $x$ if and only if $\mathcal{W}_{(0,x)} = 0$. Varying $x \in M$ completes the proof.
\end{proof}

\begin{remark}
Proposition~\ref{prop:weyl} also follows from \cite[Proposition~9.2]{InagakiJavaloyesMarcal2026} (which relates flag-plane independence of $F$ to transverse-direction independence of the null sectional curvature of $\hat{g}$) together with Lemma~\ref{lem:weyl-null-screen}, as well as from the curvature identities of Makhmali and Sagerschnig \cite{MakhmaliSagerschnig2026}. In the notation of \cite[Corollary~5.21 and equations~(5.103), (5.116)]{MakhmaliSagerschnig2026}, the invariant $\mathbf{A}$ vanishes for a Randers metric, and $\mathbf{T}$ is the trace-free part of its flag curvature; thus scalar flag curvature is equivalent to conformal flatness of the stationary lift.
\end{remark}

\section{Local normal form classification}\label{sec:normal-form}

In this section, we establish the explicit local coordinate normal form for Randers metrics of scalar flag curvature, completing the proof of Theorem~\ref{thm:main}. To deduce this normal form from Proposition~\ref{prop:weyl}, we set up three ingredients: the local Weyl--Schouten conformal-flattening criterion (cf.\ \cite[Theorem~7.37(a)]{Lee2018}), the classification of conformal Killing fields on flat space (cf.\ \cite[Chapter~1]{Schottenloher2008} and \cite[eq.~(57)]{GibbonsEtAl2009}), and unit-time normalization along spacelike sections.

\begin{lemma}[{\cite[Theorem~7.37(a)]{Lee2018}}]\label{lem:conformal-flat}
A smooth pseudo-Riemannian metric of dimension $m\ge4$ is locally conformal to a flat metric of the same signature if and only if its Weyl tensor vanishes.
\end{lemma}

We now work in Minkowski space $\mathbb{R}^{1,n}$ with coordinates $z = (\tau, x)$ and metric
\[
\eta = -d\tau^2 + \sum_{i=1}^n (dx^i)^2.
\]

\begin{lemma}[{\cite[Chapter~1]{Schottenloher2008}, \cite[eq.~(57)]{GibbonsEtAl2009}}]\label{lem:ckf}
On a connected open subset of $\mathbb{R}^{1,n}$ with $m = n + 1 \ge 3$, a smooth vector field $V$ satisfies the conformal Killing equation $\mathcal{L}_V \eta = 2\psi \eta$ for some smooth function $\psi$ if and only if $V$ has the form
\begin{equation}\label{eq:ckf}
V(z) = P_0 + Lz + \sigma z + 2\eta(q_0, z)z - \eta(z, z)q_0,
\end{equation}
where $P_0, q_0 \in \mathbb{R}^{1,n}$ are constant vectors, $\sigma \in \mathbb{R}$, and $L \in \mathfrak{so}(1,n)$ is $\eta$-skew-adjoint. In this case, the conformal factor is $\psi = \sigma + 2\eta(q_0, z)$.
\end{lemma}

Let $V$ be a timelike conformal Killing field of $\eta$ on an open subset of $\mathbb{R}^{1,n}$, and set $\Lambda = -\eta(V,V) > 0$. Along the hyperplane $\{\tau = 0\}$, decompose $V(0, x) = v(x)\partial_\tau + w(x)$ into its temporal component $v(x) \in \mathbb{R}$ and spatial component $w(x) \in \mathbb{R}^n$. Restricting $\Lambda$ to $\{\tau = 0\}$ yields
\[
D(x) := \Lambda(0, x) = -\eta\bigl(V(0, x), V(0, x)\bigr) = v(x)^2 - |w(x)|^2 > 0.
\]
Choose the time orientation determined by $V$. On each connected coordinate neighborhood, reversing the $\tau$-axis if necessary gives $v>0$, and then $D>0$ gives $v>|w|$.

Let $(s, x)$ be local flow coordinates adapted to $V$ via $(s, x) \mapsto \Phi_s(0, x)$, so that $V = \partial_s$ and $\{s = 0\} = \{\tau = 0\}$. Every spacelike hypersurface $\Sigma$ is transverse to $V$, and the future-directed null curves of $g = \eta/\Lambda$ determine a smooth, strongly convex Randers metric on $\Sigma$ via the arrival-time parametrization along the flow of $V$. On $\{s = 0\}$, set $\alpha^{(0)}(x,y) = \sqrt{a^{(0)}_x(y,y)}$, where
\begin{equation}\label{eq:a0beta0}
a^{(0)}_x(y,y) = \frac{D(x)|y|^2 + \langle w(x), y \rangle^2}{D(x)^2}, \qquad \beta^{(0)}_x(y) = \frac{\langle w(x), y \rangle}{D(x)}.
\end{equation}

The transformation of the Fermat metric under a change of spacelike section is described in \cite[Propositions~5.8 and~5.9]{CaponioJavaloyesSanchez2011}; in particular, the metric induced on the graph $s = -f(x)$ is $F_0 + df$. The following lemma expresses the spacelike condition for this graph in terms of the normal-form data $v$, $w$, and $f$.

\begin{lemma}\label{lem:unit-norm}
Under the notation above, $g = \eta/\Lambda$ admits $V$ as a unit timelike Killing field, and its future-directed null curves induce on $\{s = 0\} = \{\tau = 0\}$ the strongly convex Randers metric $F_0 = \alpha^{(0)} + \beta^{(0)}$ (the first term of \eqref{eq:F}). Moreover, for a smooth function $f$, the graph $\{s = -f(x)\}$ is spacelike if and only if \eqref{eq:convexity} holds, in which case its induced Randers metric is $F = \alpha^{(0)} + \beta^{(0)} + df$, yielding \eqref{eq:F}.
\end{lemma}
\begin{proof}
From $\mathcal{L}_V \eta = 2\psi\eta$ and $[V,V] = 0$, we have $V(\Lambda) = 2\psi\Lambda$, whence $\mathcal{L}_V(\eta/\Lambda) = 0$ and $g(V,V) = -1$. The local flow of $V$ defines adapted coordinates $(s, x)$ near $\{\tau = 0\}$ in which $V = \partial_s$ and the coefficients of $g$ are independent of $s$. Along $\{\tau = 0\}$,
\[
g_{ss} = -1, \qquad g_{si} = \frac{w_i}{D}, \qquad g_{ij} = \frac{\delta_{ij}}{D}.
\]
Completing the square yields $g = a^{(0)} - (ds - \beta^{(0)})^2$ with $a^{(0)}, \beta^{(0)}$ as in \eqref{eq:a0beta0}. The future-directed null condition $g(\dot{s}\partial_s + y, \dot{s}\partial_s + y) = 0$ has the unique positive root $\dot{s} = \sqrt{a^{(0)}(y,y)} + \beta^{(0)}(y)$, and spacelikeness of the slice is equivalent to $|\beta^{(0)}|_{a^{(0)}} < 1$. Finally, changing the time coordinate to $t = s + f(x)$ transforms the metric to
\[
g = a^{(0)} - \bigl(dt - (\beta^{(0)} + df)\bigr)^2,
\]
whose $\{t = 0\}$ slice is $s = -f(x)$. As noted after \eqref{eq:lift}, the induced metric $a^{(0)} - (\beta^{(0)} + df) \otimes (\beta^{(0)} + df)$ on this slice is positive-definite (so the slice is spacelike) if and only if $|\beta^{(0)} + df|_{a^{(0)}} < 1$. Using the inverse metric $(a^{(0)})^{ij} = D(\delta^{ij} - v^{-2}w^i w^j)$ of \eqref{eq:a0beta0} together with $D = v^2 - |w|^2 > 0$, we compute
\[
|\beta^{(0)} + df|_{a^{(0)}}^2
= \frac{|w|^2 + D\bigl(2\langle w, \nabla f \rangle + v^2 |\nabla f|^2 - \langle w, \nabla f \rangle^2\bigr)}{v^2},
\]
so $|\beta^{(0)} + df|_{a^{(0)}} < 1$ is equivalent to \eqref{eq:convexity}. Locally, subtracting a constant from $f$ ensures that the graph lies inside the flow-coordinate neighborhood.
\end{proof}

Combining Proposition~\ref{prop:weyl} with Lemmas~\ref{lem:conformal-flat}--\ref{lem:unit-norm}, we are now ready to prove the local normal form classification.

\begin{proof}[Proof of Theorem~\ref{thm:main}]
Suppose first that $F = \alpha + \beta$ has scalar flag curvature, and fix an arbitrary point $x_0 \in M$. By Proposition~\ref{prop:weyl}, the canonical Lorentzian lift $(\hat{M}, \hat{g})$ in \eqref{eq:lift} has vanishing Weyl tensor. Since $\dim \hat{M} = n + 1 \ge 4$, Lemma~\ref{lem:conformal-flat} provides an open neighborhood $\hat{U} \subset \hat{M}$ of $(0, x_0)$, a smooth function $\phi \in C^\infty(\hat{U})$, and a diffeomorphism $\Psi \colon \hat{U} \to \mathcal{O}$ onto a connected open subset $\mathcal{O} \subset \mathbb{R}^{1,n}$ such that $\Psi^* \eta = e^{2\phi} \hat{g}$, where $\eta = -d\tau^2 + \sum_{i=1}^n (dx^i)^2$ denotes the Minkowski metric. Setting the smooth positive conformal factor $\Omega := e^{-\phi \circ \Psi^{-1}}$ on $\mathcal{O}$, the metric $g := (\Psi^{-1})^* \hat{g}$ takes the form $g = \Omega^2 \eta$.

Let $V := \Psi_* T$ on $\mathcal{O}$, where $T = \partial_t$ is the unit timelike Killing field of $\hat{g}$. Since $\Psi \colon (\hat{U}, \hat{g}) \to (\mathcal{O}, g)$ is an isometry, $V$ satisfies $\mathcal{L}_V(\Omega^2 \eta) = 0$ and $(\Omega^2 \eta)(V, V) = -1$. Expanding the Lie derivative yields $\mathcal{L}_V \eta = -2V(\log \Omega)\,\eta$, so $V$ is a conformal Killing field of $\eta$. Furthermore, $\Omega^2 \eta(V, V) = -1$ shows that $V$ is timelike and $\Omega^2 = 1/\Lambda$, where $\Lambda := -\eta(V, V) > 0$; hence $g = \eta/\Lambda$. Since $\mathcal{O}$ is connected, Lemma~\ref{lem:ckf} implies that $V$ is given on $\mathcal{O}$ by \eqref{eq:ckf} for constants $P_0, q_0 \in \mathbb{R}^{1,n}$, $\sigma \in \mathbb{R}$, and $L \in \mathfrak{so}(1,n)$.

Translate $\Psi(0,x_0)$ to the origin, then apply a Lorentz transformation sending the unit timelike vector $V(0,0)/\sqrt{\Lambda(0,0)}$ to $\partial_\tau$. In these coordinates,
\[
\Psi(0, x_0) = (0, 0) \qquad \text{and} \qquad V(0, 0) = \sqrt{\Lambda(0, 0)}\,\partial_\tau.
\]
Decompose $P_0 = (a_0, p)$ and $q_0 = (-c, q)$ with $a_0, c \in \mathbb{R}$ and $p, q \in \mathbb{R}^n$. Relative to $\eta = \operatorname{diag}(-1, I_n)$, the skew-adjointness condition $\eta L + L^{\mathsf{T}} \eta = 0$ identifies $L \in \mathfrak{so}(1,n)$ uniquely with a pair $(b, Q) \in \mathbb{R}^n \times \mathfrak{so}(n)$ via
\begin{equation}\label{eq:L-decomp}
L(\tau, x) = \bigl(\langle b, x \rangle,\, b\tau + Qx\bigr).
\end{equation}
Here $P_0=V(0,0)$, so our choice of coordinates gives $p=0$ and $a_0=\sqrt{\Lambda(0,0)}>0$. We retain $p$ in the formula to include charts with an arbitrary timelike value of $V$ at the basepoint. Substituting $\eta(q_0,z)=c\tau+\langle q,x\rangle$ and $\eta(z,z)=-\tau^2+|x|^2$ in \eqref{eq:ckf} and restricting to $\tau=0$ gives
\begin{align*}
V(0, x) &= (a_0, p) + \bigl(\langle b, x \rangle,\, Qx\bigr) + \sigma(0, x) + 2\langle q, x \rangle(0, x) - |x|^2(-c, q) \\
&= \bigl(a_0 + \langle b, x \rangle + c|x|^2,\; p + Qx + (\sigma + 2\langle q, x \rangle)x - |x|^2 q\bigr) = \bigl(v(x), w(x)\bigr),
\end{align*}
recovering \eqref{eq:v}--\eqref{eq:w}. Because $v(0) = \sqrt{\Lambda(0, 0)} > 0$ and $\Lambda(0, x) = v(x)^2 - |w(x)|^2 = D(x) > 0$ near $0$, there exists an open neighborhood $U \subset \mathbb{R}^n$ of $0$ on which $v > |w|$.

Let $(s, x) \mapsto \Phi_s(0, x)$ denote the local flow coordinates adapted to $V$ based on $\{\tau = 0\}$, so that $V = \partial_s$ and $\{s = 0\} = \{\tau = 0\}$. Since $|\beta|_a < 1$, the slice $\{t = 0\} \subset \hat{M}$ is spacelike for $\hat{g}$, so its image $\Sigma_0 := \Psi\bigl(\hat{U} \cap (\{0\} \times M)\bigr)$ is a spacelike hypersurface through $(0, 0)$ in $(\mathcal{O}, g)$ and hence transverse to the timelike field $V = \partial_s$. By the implicit function theorem, shrinking $U$ around $0$ if necessary, $\Sigma_0$ is represented in flow coordinates as a graph $s = -f(x)$ for a smooth function $f \in C^\infty(U)$ with $f(0) = 0$. Since the flow-coordinate projection $(s, x) \mapsto x$ is a left inverse to $x \mapsto \Phi_{-f(x)}(0, x)$, the latter map is a diffeomorphism from $U$ onto an open neighborhood of $(0, 0)$ in $\Sigma_0$. As $\Psi$ maps $\hat{U} \cap (\{0\} \times M)$ diffeomorphically onto $\Sigma_0$, shrinking $U$ around $0$ ensures that the map $\varphi \colon U \to M$ defined by $(0, \varphi(x)) = \Psi^{-1}\bigl(\Phi_{-f(x)}(0, x)\bigr)$ is a diffeomorphism onto an open neighborhood of $x_0$ in $M$; thus $\varphi^{-1}$ defines a smooth local coordinate chart around $x_0$.
Under the coordinate change $t = s + f(x)$ in Lemma~\ref{lem:unit-norm}, the pullback of $g$ therefore agrees with the canonical lift \eqref{eq:lift} in the local coordinates $x \in U$, yielding $a = a^{(0)}$ and $\beta = \beta^{(0)} + df$. Consequently, $F$ takes the form \eqref{eq:F} on $U$, and spacelikeness of $\Sigma_0$ is equivalent to \eqref{eq:convexity}.

Conversely, suppose that constants $(a_0, c, \sigma, b, p, q, Q)$ with $Q \in \mathfrak{so}(n)$, an open set $U \subset \mathbb{R}^n$, and a function $f \in C^\infty(U)$ satisfy $v > |w|$ and \eqref{eq:convexity}. Define $P_0 = (a_0, p)$, $q_0 = (-c, q)$, and $L \in \mathfrak{so}(1,n)$ by \eqref{eq:L-decomp}, and let $V$ be the polynomial vector field \eqref{eq:ckf} on $\mathbb{R}^{1,n}$. By Lemma~\ref{lem:ckf}, $V$ is a conformal Killing field of $\eta$. Along $\{0\} \times U$, we have $V(0, x) = (v(x), w(x))$ with $v > |w| \ge 0$, so $-\eta(V, V) = v^2 - |w|^2 = D > 0$ there. By continuity, $\Lambda := -\eta(V, V)$ remains positive on an open neighborhood $\mathcal{O} \subset \mathbb{R}^{1,n}$ of $\{0\} \times U$, so $V$ is timelike on $\mathcal{O}$. Moreover, $D > 0$ ensures that $a^{(0)}$ and $\beta^{(0)}$ in \eqref{eq:a0beta0} are smooth on $U$, and the norm computation in Lemma~\ref{lem:unit-norm} shows that \eqref{eq:convexity} is equivalent to $|\beta^{(0)} + df|_{a^{(0)}} < 1$; thus \eqref{eq:F} defines a smooth, strongly convex Randers metric $F$ on $U$.

Fix $x_0\in U$ and put $\tilde f=f-f(x_0)$. Transversality of $V$ to $\{\tau=0\}$ at $(0,x_0)$ gives a flow box in which $(t,x)\mapsto\Phi_{t-\tilde f(x)}(0,x)$ is a diffeomorphism onto its image for $(t,x)\in(-\varepsilon,\varepsilon)\times U_0$, after shrinking $U_0$ about $x_0$. By Lemma~\ref{lem:unit-norm}, the pullback of $\eta/\Lambda$ is $a^{(0)}-(dt-(\beta^{(0)}+df))^2$, the canonical lift of $F$. Its Weyl tensor vanishes on this flow box by conformal flatness and on $\mathbb{R}\times U_0$ by time-translation invariance. Proposition~\ref{prop:weyl} gives scalar flag curvature on $U_0$. As $x_0$ was arbitrary, this proves the converse.
\end{proof}

\begin{remark}\label{rem:normalization-dim2}
\textup{(i)} Applying a further homothety $z \mapsto a_0^{-1}z$ of $\mathbb{R}^{1,n}$ rescales both $\eta$ and $\Lambda = -\eta(V,V)$ by $a_0^{-2}$ and hence preserves $g = \eta/\Lambda$. Combined with the translation and Lorentz boost in the proof of Theorem~\ref{thm:main}, one may therefore normalize $p = 0$ and $a_0 = 1$ in \eqref{eq:v}--\eqref{eq:w}; retaining $(a_0, p)$ accommodates charts in which $V$ takes an arbitrary timelike value at the origin.

\textup{(ii)} When $n = 2$, the canonical lift $\hat{g}$ has dimension $m = 3$, so its Weyl tensor vanishes identically and Lemma~\ref{lem:conformal-flat} no longer applies. This reflects the classical fact that every $2$-dimensional Finsler metric automatically has scalar flag curvature. Nevertheless, the sufficiency direction of Theorem~\ref{thm:main} remains valid for $n = 2$ (using the $m = 3$ case of Lemma~\ref{lem:ckf}).
\end{remark}

\section{Local Isometries and Gluing}\label{sec:global}

The local normal form in Theorem~\ref{thm:main} is not unique: in constructing the conformal Minkowski chart $\Psi$ in Section~\ref{sec:normal-form}, any post-composition of $\Psi$ with a local conformal diffeomorphism of $(\mathbb{R}^{1,n}, \eta)$ produces another valid chart, inducing a finite-dimensional gauge redundancy on the classifying data $(a_0, c, \sigma, b, p, q, Q, f)$. To describe this equivalence geometrically and pass from local charts to global manifolds, we encode each local datum in $\mathbb{R}^{1,n}$ as a pair $(V, j)$ consisting of a timelike conformal Killing field and a spacelike section.

Specifically, given a classifying datum on an open neighborhood $U \subset \mathbb{R}^n$ of a basepoint $x_0 \in U$, replace $f$ by $f - f(x_0)$ near $x_0$; this leaves $df$ and the Randers metric \eqref{eq:F} unchanged while normalizing $f(x_0) = 0$. Let $V$ be the polynomial conformal Killing field \eqref{eq:ckf} determined by $(a_0, c, \sigma, b, p, q, Q)$ via \eqref{eq:L-decomp}, with local flow $\Phi_s$. Since $V$ is timelike near $(0, x_0)$ and $f(x_0) = 0$, the spacelike section embedding
\begin{equation}\label{eq:section-j}
j(x) := \Phi_{-f(x)}(0, x)
\end{equation}
is defined on a sufficiently small neighborhood of $x_0$ and satisfies $j(x_0) = (0, x_0)$. Together with $\Phi$, it defines the adapted local diffeomorphism $\Theta(t, x) := \Phi_t(j(x)) = \Phi_{t - f(x)}(0, x)$ onto an open subset of $\mathbb{R}^{1,n}$ pulling back the unit-normalized metric $\eta/\Lambda$ to the canonical Lorentzian lift \eqref{eq:lift}, sending $\partial_t$ to $V$, and identifying the slice $\{t = 0\}$ with the image of $j$. Note that $\Phi$ is uniquely determined by an explicit polynomial ODE and imposes no additional constraint on the classifying data.

Javaloyes, Lichtenfelz, and Piccione \cite{JavaloyesLichtenfelzPiccione2015} related conformal maps between stationary spacetimes to almost isometries of their Fermat metrics. Let $\Psi$ be a conformal map carrying the chosen stationary field $T$ to the chosen stationary field $T'$. In coordinates adapted to these fields, $\Psi$ takes the form $(t,x)\mapsto(t+u(x),\varphi(x))$. If $\Psi$ also preserves the marked zero sections, then $u=0$ and $\varphi$ is an isometry of the induced Randers metrics. We express this condition for the pairs $(V,j)$ via local conformal transformations of Minkowski space.

\begin{lemma}[Liouville's theorem in Minkowski space; {\cite[Theorem~1.9]{Schottenloher2008}}]\label{lem:liouville}
Let $m = n + 1 \ge 3$ and let $\mathcal{O} \subset \mathbb{R}^{1,n}$ be a connected open subset. Every smooth conformal diffeomorphism $\Psi \colon \mathcal{O} \to \Psi(\mathcal{O}) \subset \mathbb{R}^{1,n}$ (satisfying $\Psi^*\eta = e^{2\omega}\eta$ for some $\omega \in C^\infty(\mathcal{O})$) takes the form
\[
\Psi(z) = a_2 + \lambda A\, K_B(z - a_1),
\]
where $a_1, a_2, B \in \mathbb{R}^{1,n}$, $\lambda > 0$, $A \in O(1,n)$, and $K_B$ is the special conformal transformation
\begin{equation}\label{eq:KB}
K_B(w) = \frac{w - \eta(w,w) B}{1 - 2\eta(B,w) + \eta(B,B)\eta(w,w)},
\end{equation}
whose denominator is nowhere vanishing for $w = z - a_1$ with $z \in \mathcal{O}$.
\end{lemma}

\begin{theorem}\label{thm:local-iso}
Let $F$ on $U \subset \mathbb{R}^n$ and $F'$ on $U' \subset \mathbb{R}^n$ be two Randers metrics arising from classifying data in Theorem~\ref{thm:main}, with associated pairs $(V, j)$ near $x_0 \in U$ and $(V', j')$ near $x'_0 \in U'$. A local diffeomorphism $\varphi \colon U_0 \to U'_0$ between open neighborhoods of $x_0$ and $x'_0$ with $\varphi(x_0) = x'_0$ represents a germ of a Randers isometry at $x_0$
(that is, $\varphi^*F'=F$ on some neighborhood of $x_0$) if and only if there exists a conformal diffeomorphism $\Psi$ of $(\mathbb{R}^{1,n}, \eta)$ defined on an open neighborhood of $j(x_0)$ with $\Psi(j(x_0)) = j'(x'_0)$ such that (after shrinking $U_0$ if necessary)
\begin{equation}\label{eq:equiv}
\Psi_* V = V', \qquad \Psi \circ j = j' \circ \varphi.
\end{equation}
In particular, any such conformal map $\Psi$ mapping $j(U_0)$ into the image of $j'$ uniquely determines the local Randers isometry via $\varphi = (j')^{-1} \circ \Psi \circ j$.
\end{theorem}

\begin{proof}
Let $F = \alpha + \beta$ and $F' = \alpha' + \beta'$ be the Randers metrics defined by the two data sets near basepoints $x_0 \in U$ and $x'_0 \in U'$, with $f(x_0) = 0$ and $f'(x'_0) = 0$. Denote by $\Theta(t, x) = \Phi_t(j(x))$ and $\Theta'(t', x') = \Phi'_{t'}(j'(x'))$ their adapted local diffeomorphisms onto open neighborhoods $\mathcal{O}, \mathcal{O}' \subset \mathbb{R}^{1,n}$ of $j(x_0)$ and $j'(x'_0)$. By Lemma~\ref{lem:unit-norm} and the proof of Theorem~\ref{thm:main},
\[
\Theta^*(\eta/\Lambda) = \hat{g}, \qquad \Theta_* \partial_t = V, \qquad \Theta(0, x) = j(x),
\]
where $\hat{g} = a - (dt - \beta)^2$ is the canonical Lorentzian lift \eqref{eq:lift} of $F$, and analogously for $(\Theta', \Lambda', V', j', \hat{g}')$.

Suppose first that $\varphi \colon U_0 \to U'_0$ is a local Randers isometry ($\varphi^* F' = F$) between neighborhoods of $x_0$ and $x'_0$ with $\varphi(x_0) = x'_0$. Because $\alpha$ and $\beta$ are uniquely recovered as the even and odd parts of $F$ in $y$,
\[
\alpha(x, y) = \frac{1}{2}\bigl(F(x, y) + F(x, -y)\bigr), \qquad \beta_x(y) = \frac{1}{2}\bigl(F(x, y) - F(x, -y)\bigr),
\]
the map $\varphi$ preserves $\alpha$ and $\beta$ separately; hence $\varphi^* a' = a$ (by polarization) and $\varphi^* \beta' = \beta$. Consequently, the product lift $\hat{\varphi}(t, x) := (t, \varphi(x))$ satisfies $\hat{\varphi}^* \hat{g}' = \hat{g}$, $\hat{\varphi}_* \partial_t = \partial_{t'}$, and $\hat{\varphi}(0, x) = (0, \varphi(x))$. Shrinking $U_0$ and the $t$-interval if necessary, the conjugated diffeomorphism $\Psi := \Theta' \circ \hat{\varphi} \circ \Theta^{-1}$ is defined near $j(x_0)$ and satisfies
\[
\Psi^*\!\left(\frac{\eta}{\Lambda'}\right) = (\Theta^{-1})^* \hat{\varphi}^* (\Theta')^*\!\left(\frac{\eta}{\Lambda'}\right) = (\Theta^{-1})^* \hat{\varphi}^* \hat{g}' = (\Theta^{-1})^* \hat{g} = \frac{\eta}{\Lambda}.
\]
Hence $\Psi^*\eta = \frac{\Lambda' \circ \Psi}{\Lambda}\eta$, showing that $\Psi$ is a local conformal diffeomorphism of $(\mathbb{R}^{1,n}, \eta)$ and therefore has the explicit form of Lemma~\ref{lem:liouville} near $j(x_0)$. Moreover,
\[
\Psi_* V = \Theta'_* \hat{\varphi}_* \Theta_*^{-1} V = \Theta'_* \hat{\varphi}_* \partial_t = \Theta'_* \partial_{t'} = V',
\]
and for every $x \in U_0$,
\[
(\Psi \circ j)(x) = (\Theta' \circ \hat{\varphi} \circ \Theta^{-1})\bigl(\Theta(0, x)\bigr) = \Theta'\bigl(0, \varphi(x)\bigr) = (j' \circ \varphi)(x),
\]
establishing \eqref{eq:equiv}.

Conversely, suppose \eqref{eq:equiv} holds for a local conformal map $\Psi$ with $\Psi^*\eta = e^{2\omega}\eta$ and a local diffeomorphism $\varphi \colon U_0 \to U'_0$. At every point $z$ in the domain of $\Psi$, the relation $\Psi_* V = V'$ yields
\begin{align*}
\Lambda'\bigl(\Psi(z)\bigr) &= -\eta\bigl(V'(\Psi(z)), V'(\Psi(z))\bigr) = -\eta\bigl(d\Psi_z(V(z)), d\Psi_z(V(z))\bigr) \\
&= -(\Psi^*\eta)_z\bigl(V(z), V(z)\bigr) = -e^{2\omega(z)}\eta\bigl(V(z), V(z)\bigr) = e^{2\omega(z)}\Lambda(z).
\end{align*}
Dividing $\Psi^*\eta = e^{2\omega}\eta$ by $\Lambda' \circ \Psi = e^{2\omega}\Lambda$ cancels the conformal factor $e^{2\omega}$, giving
\[
\Psi^*\!\left(\frac{\eta}{\Lambda'}\right) = \frac{\eta}{\Lambda}.
\]
Furthermore, since $\Psi_* V = V'$, $\Psi$ intertwines the local flows $\Phi_t$ and $\Phi'_t$ ($\Psi \circ \Phi_t = \Phi'_t \circ \Psi$ for small $|t|$). Combined with $\Psi \circ j = j' \circ \varphi$, this implies
\[
\Psi\bigl(\Theta(t, x)\bigr) = \Psi\bigl(\Phi_t(j(x))\bigr) = \Phi'_t\bigl(\Psi(j(x))\bigr) = \Phi'_t\bigl(j'(\varphi(x))\bigr) = \Theta'\bigl(t, \varphi(x)\bigr) = (\Theta' \circ \hat{\varphi})(t, x)
\]
for $(t, x)$ near $(0, x_0)$. Pulling back $\eta/\Lambda'$ along both sides gives
\[
\hat{\varphi}^* \hat{g}' = \hat{\varphi}^* (\Theta')^*\!\left(\frac{\eta}{\Lambda'}\right) = \Theta^* \Psi^*\!\left(\frac{\eta}{\Lambda'}\right) = \Theta^*\!\left(\frac{\eta}{\Lambda}\right) = \hat{g}.
\]
Since $\hat{\varphi}(t, x) = (t, \varphi(x))$ preserves $dt$, comparing the cross terms $dt \otimes dx^i$ and spatial terms $dx^i \otimes dx^j$ in $\hat{\varphi}^* \bigl(a' - (dt - \beta')^2\bigr) = a - (dt - \beta)^2$ yields $\varphi^* \beta' = \beta$ and $\varphi^* a' = a$ (equivalently, $\Psi$ preserves future-directed null curves and arrival-time increments in Lemma~\ref{lem:unit-norm}). Hence $\varphi^* F' = F$, so $\varphi$ is a local Randers isometry. Finally, because $j'$ is an embedding onto its image, $\Psi \circ j = j' \circ \varphi$ uniquely forces $\varphi = (j')^{-1} \circ \Psi \circ j$.
\end{proof}

To pass from the local parametrization to a global description, define a \emph{normal-form atlas} on a smooth $n$-manifold $M$ ($n \ge 3$) to be a smooth atlas $\{(V_\alpha, \chi_\alpha)\}_{\alpha \in I}$ such that each coordinate domain $U_\alpha := \chi_\alpha(V_\alpha) \subset \mathbb{R}^n$ carries classifying data $(a_{0,\alpha}, c_\alpha, \sigma_\alpha, b_\alpha, p_\alpha, q_\alpha, Q_\alpha, f_\alpha)$ from Theorem~\ref{thm:main}, and each transition diffeomorphism $\varphi_{\alpha\gamma} := \chi_\gamma \circ \chi_\alpha^{-1}$ locally satisfies \eqref{eq:equiv} near every point $x_0 \in \chi_\alpha(V_\alpha \cap V_\gamma)$ for an explicit conformal map $\Psi_{\alpha\gamma}$ from Lemma~\ref{lem:liouville}, where the local section maps $j$ and $j'$ in \eqref{eq:equiv} are normalized at $x_0$ and $x'_0 := \varphi_{\alpha\gamma}(x_0)$ according to \eqref{eq:section-j}:
\[
j_{\alpha, x_0}(x) := \Phi^\alpha_{-(f_\alpha(x) - f_\alpha(x_0))}(0, x), \qquad
j_{\gamma, x'_0}(x') := \Phi^\gamma_{-(f_\gamma(x') - f_\gamma(x'_0))}(0, x').
\]
Equivalently, such a manifold arises abstractly by gluing open sets $\{U_\alpha \subset \mathbb{R}^n\}_{\alpha \in I}$ with data from Theorem~\ref{thm:main} along cocycle-compatible transition maps locally satisfying \eqref{eq:equiv}, provided the quotient space is Hausdorff and second-countable.

Because strong convexity and scalar flag curvature are local properties and Theorem~\ref{thm:local-iso} identifies maps locally obeying \eqref{eq:equiv} with local Randers isometries, the local metrics \eqref{eq:F} agree on all overlaps and patch to a smooth global metric; conversely, covering a given manifold by coordinate charts from Theorem~\ref{thm:main} yields a normal-form atlas. Similarly, any family of local diffeomorphisms satisfying \eqref{eq:equiv} (with basepoint-normalized sections as above) on refined atlases that commutes with the transition maps glues to a smooth local isometry, which is a global isometry once it is bijective (by the inverse function theorem). Theorems~\ref{thm:main} and~\ref{thm:local-iso} therefore yield the following immediate global description.

\begin{corollary}[Global gluing description]\label{cor:global-atlas}
Let $n \ge 3$.
\begin{enumerate}
\item[\textup{(i)}] On any smooth $n$-manifold $M$ equipped with a normal-form atlas, the local formulas \eqref{eq:F} patch together to define a smooth, strongly convex Randers metric $F$ of scalar flag curvature. Conversely, every smooth, strongly convex Randers $n$-manifold $(M, F)$ of scalar flag curvature admits a normal-form atlas inducing $F$.
\item[\textup{(ii)}] Two such manifolds $(M, F)$ and $(M', F')$ are globally isometric if and only if, after refining their normal-form atlases, there exists a family of local diffeomorphisms satisfying \eqref{eq:equiv} that commutes with the transition maps and induces a bijection $M \to M'$.
\end{enumerate}
\end{corollary}

\begin{remark}\label{rem:global-cohomology}
Corollary~\ref{cor:global-atlas} ensures that the local Randers metrics patch together to define $F$. If the local background metrics $F_{0,\alpha}$ also agree on overlaps, they assemble into a global metric $F_0$ (of scalar flag curvature and isotropic $S$-curvature), and the local exact $1$-forms $df_\alpha = F - F_{0,\alpha}$ patch together into a global closed $1$-form $\theta$ satisfying $F = F_0 + \theta$. For such a choice of $F_0$, a global potential $f \in C^\infty(M)$ with $\theta = df$ exists if and only if $[\theta] = 0$ in $H^1_{\mathrm{dR}}(M)$. In particular, when $H^1_{\mathrm{dR}}(M) = 0$, every compatible family $\{F_{0,\alpha}\}$ yields a global decomposition $F = F_0 + df$.
\end{remark}

\section{Examples and Navigation Interpretation}\label{sec:examples}

In this section, we illustrate Theorem~\ref{thm:main} with two explicit examples exhibiting direction-dependent flag curvature and non-isotropic $S$-curvature, derive the Zermelo navigation data of the general normal form, and prove Corollary~\ref{cor:projectively-isotropic}.

\begin{example}[Non-constant scalar flag curvature]\label{ex:nonconst}
In Theorem~\ref{thm:main}, take $v = 1$, $w = 0$, and $f(x) = \frac{\epsilon}{2}(x^1)^2$ with $\epsilon \neq 0$. On the strip $|\epsilon x^1| < 1$, the strong convexity condition \eqref{eq:convexity} holds and \eqref{eq:F} gives
\[
F(x,y) = |y| + \epsilon x^1 y^1.
\]
Under the time change $s = t - f(x)$, its canonical Lorentzian lift \eqref{eq:lift} is the flat Minkowski metric $|dx|^2 - ds^2$. Since $d\beta = 0$ and $a$ is Euclidean, the magnetic spray coefficients in \eqref{eq:mag-spray} vanish in Cartesian coordinates, so the geodesic spray of $F$ is purely projective: $G^i = P y^i$ for a positively $1$-homogeneous function $P(x,y)$. Along an $F$-geodesic with constant $F$-speed, differentiating $F(x,\dot{x})$ gives $0 = F_{x^k} y^k - 2PF$, hence
\[
P(x,y) = \frac{\epsilon (y^1)^2}{2 F(x,y)}.
\]
The magnetic spray is $\mathcal{S}_{\mathrm{mag}}= y^k\partial_{x^k}$, so its Jacobi endomorphism vanishes and $\mathcal{S}_{\mathrm{mag}}P=y^kP_{x^k}$. Hence Lemma~\ref{lem:proj-invariance} gives $R_y \equiv (P^2-y^kP_{x^k})\,\Id \pmod{\mathbb{R}y}$. Using $F_{x^1} = \epsilon y^1$, we compute
\[
P^2 - y^k P_{x^k}
= \frac{\epsilon^2 (y^1)^4}{4 F^2} + \frac{\epsilon (y^1)^2}{2 F^2} (\epsilon y^1) y^1
= \frac{3\epsilon^2 (y^1)^4}{4 F^2},
\]
and dividing by $F(x,y)^2$ gives the scalar flag curvature
\[
K(x,y) = \frac{3\epsilon^2 (y^1)^4}{4 F(x,y)^4}.
\]
In particular, at the origin we have $K(0, e_1) = 3\epsilon^2/4$ whereas $K(0, e_2) = 0$. Thus a flat Lorentzian lift admits a spacelike slice whose induced Randers metric has genuinely direction-dependent scalar flag curvature. This example also follows from \cite[eq.~(1.11)]{ShenXia2014} by taking $\alpha$ to be Euclidean and $f(x)=\epsilon(x^1)^2/2$.
\end{example}

\begin{example}[Non-projectively flat metric with non-isotropic $S$-curvature, \cite{ChenZhao2012}]\label{ex:chen-zhao}
Let $n = 2N \ge 4$, identify $\mathbb{R}^{2N} \cong \mathbb{C}^N$ with its standard complex structure $J \in \mathfrak{so}(2N)$ ($J^2 = -\Id$), and choose in Theorem~\ref{thm:main} the data
\[
v(x) = 1 + \frac{1}{4}|x|^2, \qquad w(x) = Jx, \qquad f(x) = \epsilon \log\!\left(\frac{1 + \frac{1}{4}|x|^2}{1 - \frac{1}{4}|x|^2}\right) \qquad (\epsilon \neq 0)
\]
on the background ball $\{|x| < 2\}$. Here $D(x) = \bigl(1 - \frac{1}{4}|x|^2\bigr)^2 > 0$, so $v > |w|$ throughout $\{|x| < 2\}$.

The radial diffeomorphism $z = x/v(x)$ maps $\{|x| < 2\}$ onto the unit ball $\mathbb{B}^{2N} = \{|z| < 1\} \subset \mathbb{C}^N$ and satisfies $1 - |z|^2 = D/v^2$. Differentiating $z$ along a tangent vector $y$ gives $dz = v^{-1}y - \frac{1}{2}v^{-2}\langle x,y\rangle x$, whence
\[
|dz|^2 = \frac{|y|^2}{v^2} - \frac{\langle x, y \rangle^2}{v^4}, \qquad
\langle z, dz \rangle = \frac{\sqrt{D}}{v^3}\langle x, y \rangle, \qquad
\langle Jz, dz \rangle = \frac{\langle Jx, y \rangle}{v^2}.
\]
Because $(1 - |z|^2)|dz|^2 + \langle z, dz \rangle^2 = (D/v^4)|y|^2$, the ingredients $a^{(0)}$ and $\beta^{(0)}$ in \eqref{eq:a0beta0} take the complex-analytic form
\[
a^{(0)}_x(y, y) = \frac{(1 - |z|^2)|dz|^2 + \langle z, dz \rangle^2 + \langle Jz, dz \rangle^2}{(1 - |z|^2)^2},
\qquad
\beta^{(0)}_x(y) = \frac{\langle Jz, dz \rangle}{1 - |z|^2}.
\]
Thus $a^{(0)}$ is the Bergman metric on $\mathbb{B}^{2N}$ of constant holomorphic sectional curvature $-4$. Setting $u(z) := -\frac{1}{2}\log(1 - |z|^2)$, we have $f(x) = \epsilon\, u(z)$ and $\beta^{(0)} = du(-J\,\cdot)$, whose exterior derivative is twice the K\"ahler form (so $d\beta \neq 0$). Since $\nabla f(x) = \epsilon x/(v\sqrt{D})$ is orthogonal to $w(x) = Jx$, the strong convexity condition \eqref{eq:convexity} simplifies to
\[
v^2 |\nabla f|^2 = \frac{\epsilon^2 |x|^2}{D} = \frac{\epsilon^2 |z|^2}{1 - |z|^2} < 1, \qquad \text{i.e.,} \qquad |z| < \frac{1}{\sqrt{1 + \epsilon^2}},
\]
which defines the domain of the resulting Randers metric $F$.

In the complex coordinates on this subball, $F(z,\xi) = \sqrt{a^{(0)}(\xi,\xi)} + du(\epsilon\xi - J\xi)$. Since complex conjugation $C(z) = \bar{z}$ preserves both $a^{(0)}$ and $u$ and satisfies $dC \circ J = -J \circ dC$, we have $F = C^* F_{\mathrm{CZ}, -\epsilon}$, where $F_{\mathrm{CZ}, -\epsilon}$ is the metric of \cite[Theorem~1.1 and Proposition~3.2]{ChenZhao2012} with potential $f_{\mathrm{CZ}} = \frac{1}{2}\log(1 - |z|^2) = -u$. By \cite{ChenZhao2012}, for $\epsilon \neq 0$ this metric is neither locally projectively flat nor of isotropic $S$-curvature.

In the notation of \cite[Proposition~3.2]{ChenZhao2012}, subtracting the exact $1$-form $\lambda\,df_{\mathrm{CZ}}$ from $F_{\mathrm{CZ},\lambda}$ yields $F_{\mathrm{CZ},0}$ throughout the subball $|z| < (1 + \lambda^2)^{-1/2}$. Since $F_{\mathrm{CZ},0}$ has scalar flag curvature and vanishing $S$-curvature \cite[pp.~481--482]{ChenZhao2012}, its flag curvature is weakly isotropic by \cite[eq.~(1.6)]{ChengShen2009}. Consequently, the Chen--Zhao metrics in this family are in fact projectively isotropic on their entire domains of strong convexity, in contrast to the remark of Shen and Xia \cite{ShenXia2013}.

Finally, we verify \eqref{eq:Qu} on this example. Writing $a = a^{(0)}$ and $\omega(X,Y) = a(JX,Y)$, the relation $d\beta = 2\omega$ gives $A = J$ and $\nabla^a A = 0$. Using the curvature tensor of a complex space form of constant holomorphic sectional curvature $-4$,
\[
R^a(X,U)U = -X - 3a(X,JU)JU \qquad \bigl(|U|_a = 1,\; X \perp_a U\bigr),
\]
together with $-A^2 X = -J^2 X = X$, we see that all terms in \eqref{eq:Qu} cancel and $Q_U X = 0$ for every $X \perp_a U$.
\end{example}

\begin{remark}[The examples of Shen and Xia]\label{rem:shen-xia-examples}
In each case below, we exhibit a representative with constant flag curvature and vanishing $S$-curvature, thereby showing that the metrics in \cite[Examples~1.1 and~1.2]{ShenXia2013} are projectively isotropic on their respective domains of strong convexity. This corrects the assertions about projective isotropy following equations~(1.19) and~(1.24) of that paper.

\smallskip
\noindent\emph{(i) The rotational example.}
In \cite[Example~1.1]{ShenXia2013}, write the wind in column-vector notation as $V(x) = Qx$ with $Q^{\mathsf T} = -Q$, and set $f(x) = \frac{1}{2}|x|^2$. On the domain $U = \{x \in \mathbb{R}^n : |x| < 1,\ |Qx| < 1\}$, let $H$ be the navigation metric on the Euclidean background with wind $V$; since $V$ is a Euclidean Killing field, $H$ has zero flag curvature and vanishing $S$-curvature \cite{BaoRoblesShen2004, ChengShen2009}.

The metric $\widetilde F$ in \cite[Example~1.1]{ShenXia2013} is defined via Zermelo navigation on the Randers metric $|y| + df_x(y)$ with the same wind $V$. Since $df_x(V(x)) = \langle x, Qx \rangle = 0$, the navigation relation for $\tau = \widetilde F(x,y)$ simplifies to
\[
  |y - \tau V(x)| = \tau - df_x(y).
\]
To absorb the exact term $df_x(y)$ into a coordinate change, consider the diffeomorphism $\Phi(x) = e^{-f(x)Q}x$ of $U$ (with inverse $x \mapsto e^{f(x)Q}x$), which preserves both $|x|$ and $|Qx|$. Setting $s = \tau - df_x(y)$ and differentiating $\Phi$ along $y$, the identities $df_x(y)Q\,\Phi(x) = df_x(y)V(\Phi(x))$ and $e^{-f(x)Q}V(x) = V(\Phi(x))$ yield
\[
  d\Phi_x(y) - sV(\Phi(x))
  = e^{-f(x)Q}\bigl(y - \tau V(x)\bigr).
\]
Taking Euclidean norms gives $|d\Phi_x(y) - sV(\Phi(x))| = s$, which is precisely the navigation equation for $s = H(\Phi(x), d\Phi_x(y))$. Hence $\tau = s + df_x(y)$ translates to
\[
  \widetilde F = \Phi^*H + df \qquad \text{on } U.
\]
Since $\Phi^*H$ is isometric to $H$, it has zero flag curvature and vanishing $S$-curvature.

\smallskip
\noindent\emph{(ii) The radial example.}
In \cite[Example~1.2]{ShenXia2013}, set $r = |x|$ and $L = (1 + 3cr)(1 - cr)$ with $c \neq 0$. On the domain $\{r > 0,\ L > 0\}$, write $\widetilde F = \sqrt{\widetilde a} + \widetilde\beta$. In polar coordinates, equations~(1.22)--(1.23) of \cite{ShenXia2013} read
\[
  \widetilde a
  = \frac{dr^2}{L^2} + \frac{r^2}{L}\,g_{\mathbb S^{n-1}},
  \qquad
  \widetilde\beta
  = \frac{1 - 3cr}{2L}\,dr
  = d\!\left(\frac{\log L}{4c}\right),
\]
where $g_{\mathbb S^{n-1}}$ is the unit round metric. Thus $\widetilde\beta$ is already exact on the entire domain. Passing to the radial coordinate $\rho = r/\sqrt{L}$ and noting that $L > 0$ forces $-\frac{1}{3} < cr < 1$ (hence $1 + cr > 0$), we compute
\[
  \frac{d\rho}{dr} = \frac{1 + cr}{L^{3/2}} > 0
  \qquad \text{and} \qquad
  1 + 4c^2\rho^2 = \frac{(1 + cr)^2}{L},
\]
which transforms $\widetilde a$ into
\[
  \widetilde a
  = \frac{d\rho^2}{1 + 4c^2\rho^2} + \rho^2 g_{\mathbb S^{n-1}}.
\]
This is the standard hyperbolic metric of constant sectional curvature $-4c^2$. Consequently, the Riemannian metric $\sqrt{\widetilde a}$ itself serves as the required constant-curvature, zero-$S$-curvature representative, giving
\[
  \widetilde F
  = \sqrt{\widetilde a} + d\!\left(\frac{\log L}{4c}\right).
\]
\end{remark}

\paragraph{Zermelo navigation interpretation.}
Recall from \cite{BaoRoblesShen2004, ChengShen2009} that a strongly convex Randers metric $F = \alpha + \beta$ (with $\alpha = \sqrt{a}$) corresponds to Zermelo navigation data $(h, W)$ via
\[
h = \varepsilon\,(a - \beta \otimes \beta), \qquad W = -\frac{\beta^\sharp}{\varepsilon},
\]
where $\varepsilon := 1 - |\beta|_a^2 > 0$ and $\beta^\sharp$ denotes the $a$-metric dual of $\beta$.

To express $(h, W)$ explicitly for the normal form \eqref{eq:F} in Theorem~\ref{thm:main}, write $a = a^{(0)}$ and $\beta = \beta^{(0)} + df$ as in \eqref{eq:a0beta0}, and set
\[
E := 1 - \bigl(2\langle w, \nabla f \rangle + v^2 |\nabla f|^2 - \langle w, \nabla f \rangle^2\bigr).
\]
By \eqref{eq:convexity}, we have $E > 0$. Using the inverse metric $(a^{(0)})^{ij} = D(\delta^{ij} - v^{-2}w^i w^j)$ of \eqref{eq:a0beta0} together with $D = v^2 - |w|^2 > 0$ and the norm identity in the proof of Lemma~\ref{lem:unit-norm}, we obtain
\[
\varepsilon = 1 - |\beta|_a^2 = \frac{D}{v^2}\, E, \qquad
\beta^\sharp = D\,\nabla f + \frac{D}{v^2}\bigl(1 - \langle w, \nabla f \rangle\bigr) w.
\]
Moreover, setting $\delta := \sum_{i=1}^n (dx^i)^2$ and $w^\flat := \langle w, \cdot \rangle$, the rank-one terms $D^{-2} w^\flat \otimes w^\flat$ cancel in $a - \beta \otimes \beta$, leaving
\[
a - \beta \otimes \beta = \frac{1}{D}\Bigl(\delta - \bigl(w^\flat \otimes df + df \otimes w^\flat\bigr) - D\, df \otimes df\Bigr).
\]
Substituting these expressions into the definitions of $h$ and $W$ yields the explicit Zermelo navigation data of $F$:
\begin{align}
h &= \frac{E}{v^2}\Bigl(\delta - \bigl(w^\flat \otimes df + df \otimes w^\flat\bigr) - D\, df \otimes df\Bigr), \label{eq:nav-h} \\
W &= -\frac{\bigl(1 - \langle w, \nabla f \rangle\bigr) w + v^2 \nabla f}{E}. \label{eq:nav-W}
\end{align}
On the time-oriented domain $v > 0$, the condition $D > 0$ is equivalent to $v > |w|$. Under these equivalent conditions, $|W|_h^2 = |\beta|_a^2 = 1 - (D/v^2)E$, so the mild-wind condition $|W|_h < 1$ is equivalent to $E > 0$.

Specializing \eqref{eq:nav-h}--\eqref{eq:nav-W} to the case $f = 0$ immediately yields the Zermelo navigation characterization of the unperturbed metric $F_0$ and proves Corollary~\ref{cor:projectively-isotropic}.

\begin{proof}[Proof of Corollary~\ref{cor:projectively-isotropic}]
Suppose first that $F$ has scalar flag curvature, and work in a local coordinate neighborhood $U$ given by Theorem~\ref{thm:main}. Let $F_0 = \alpha^{(0)} + \beta^{(0)}$ be the metric obtained by setting $f = 0$ in \eqref{eq:F}. Since $v > |w|$, $F_0$ is a smooth, strongly convex Randers metric satisfying $F = F_0 + df$. Setting $f = 0$ in \eqref{eq:nav-h}--\eqref{eq:nav-W} (so that $E = 1$) gives the Zermelo navigation data of $F_0$:
\[
h_0 = v^{-2}\delta, \qquad W_0 = -w.
\]
In Euclidean coordinates, $v(x) = a_0 + \langle b, x \rangle + c|x|^2$ satisfies $\nabla^2 v = 2c\,\delta$ and $4cv - |\nabla v|^2 = 4a_0 c - |b|^2 =: \kappa$. By the conformal transformation formula for the Schouten tensor \cite[Theorem~7.30]{Lee2018}, the Schouten tensor of $h_0 = v^{-2}\delta$ is
\[
\mathsf{P}^{h_0} = v^{-1}\nabla^2 v - \frac{1}{2}v^{-2}|\nabla v|^2\delta = \frac{\kappa}{2}\,h_0.
\]
Since $h_0$ is conformally flat, it has constant sectional curvature $\kappa$. Moreover, differentiating \eqref{eq:w} yields $\mathcal{L}_w \delta = 2(\sigma + 2\langle q, x \rangle)\delta$, so $W_0 = -w$ is conformal with respect to $h_0$:
\[
\mathcal{L}_{W_0} h_0 = 2\rho\, h_0, \qquad \rho = -\sigma - 2\langle q, x \rangle + \frac{w(v)}{v}.
\]
By the characterization of Cheng and Shen \cite[Theorem~1.1]{ChengShen2009}, $F_0$ has scalar flag curvature and isotropic $S$-curvature, specifically $S_0 = (n+1)c_S F_0$ with $c_S = -\rho/2$.

Conversely, suppose that every point admits a neighborhood on which $F = F_0 + df$, where $f$ is smooth and $F_0$ is a strongly convex Randers metric of scalar flag curvature and isotropic $S$-curvature. Since their length functionals differ only by a boundary term, $F$ and $F_0$ have the same oriented unparametrized geodesics, meaning that their geodesic sprays are projectively related; Lemma~\ref{lem:proj-invariance} therefore implies that $F$ has scalar flag curvature. By \cite[Theorem~1.1]{ChenMoShen2003} and \cite[eq.~(1.6)]{ChengShen2009}, the identity
$S_0=(n+1)c_S F_0$ gives
\[
  K_0(x,y)=\frac{3(dc_S)_x(y)}{F_0(x,y)}+\sigma_0(x)
\]
for some smooth function $\sigma_0$, so $F_0$ has weakly isotropic flag curvature. Thus $F=F_0+df$ is locally projectively isotropic in the sense of \cite{ShenXia2013}.
\end{proof}

\vskip 2mm
\paragraph{Declaration of generative AI.} The author used generative AI as an assisting tool in the development of this manuscript and takes full responsibility for the content of the work.


\vskip 2mm
\noindent
Yalu Feng \\
School of Mathematical Sciences, Xiamen University, Xiamen 361005, P.R. China \\
E-mail: \texttt{fengyl2824@qq.com}

\end{document}